\documentclass{article}
\usepackage[utf8]{inputenc}
\usepackage{mathtools} 
\usepackage{amsfonts} 
\usepackage{amsmath}
\usepackage{float}
\usepackage{nccmath}
\usepackage[english]{babel}

\title{On the $p$-adic valuation of a hyperfactorial}
\author{Luca Onnis}
\begin{document}
\date{November 2020}

\maketitle
\begin{abstract}
    In this document will be proved a formula to compute the p-adic valuation of a hyperfactorial. We call a hyperfactorial the result of multiplying a given number of consecutive integers from 1 to the given number,each raised to its own power. For example, the hyperfactorial of n is equal to: $1^1 2^2 3^3\dots n^n$ . Lots of studies have been done about the hyperfactorial function, in particular two mathematicians: Glaisher and Kinkelin, who have found the asymptotic behaviour of this function as n that approaches infinity (finding a costant, the Glaisher-Kinkelin costant \cite{cg}, which has a lot of expressions using the Euler Gamma function and the Riemann Zeta function). In particular in this document I'll write about the p-adic valuation of this function, or rather the maximum exponent of $p$($p$ a prime integer) such that $p$ raised to that power divides the hyperfactorial of $n$. The formula which I will present uses the famous De-Polignac formula for the p-adic valuation of the simple factorial. Then I'll discuss about the asymptotic analysis of our result.
\end{abstract}

\section{Introduction}
I define $H_f(n)$ as the hyperfactorial function, where:
$H_f(n)=1^1 2^2 3^3\dots n^n$
Let's assume I want to compute the 5-adic valuation of the hyperfactorial of 60. With my formula, I can say that $5^{465}\mid H_f(60)$ while $5^{466}\not | H_f(60)$.
$a\mid b$ indicate that $b$ is divisible by $a$ \\
Let's consider the 2-adic valuation of the hyperfactorial of 10, this function satisfies an important property,in fact we have that: $v_p(ab)=v_p(a)+v_p(b)$
Hence, in general:
\[
v_p[H_f(n)]=v_p(1^1)+v_p(2^2)+v_p(3^3)+\dots+v_p(n^n)
\]
But applying again the same property we'll have:
\[
v_p[H_f(n)]=1v_p(1)+2v_p(2)+3v_p(3)+\dots+nv_p(n)
\]
And:
\[
v_2[H_f(10)]=1v_2(1)+2v_2(2)+3v_3(3)+\dots+10v_2(10)
\]
It's evident how the 2-adic valuation of odd numbers is 0, so we can eliminate those terms \\
We want to re-write it to have only simple factorial as the arguments of the p-adic valuation, so:
\[
v_2(10!)=v_2(1)+v_2(2)+\dots+v_2(10)
\]
Hence:
\[
v_2[H_f(10)]=10v_2(10!)-2v_2(8)-4v_2(6)-6v_2(4)-8v_2(2)
\]
\[
v_2[H_f(10)]=10v_2(10!)-2v_2(8!)-2v_2(6)-4v_2(4)-6v_2(2)
\]
\[
v_2[H_f(10)]=10v_2(10!)-2v_2(8!)-2v_2(6!)-2v_2(4)-4v_2(2)
\]
\[
v_2[H_f(10)]=10v_2(10!)-2v_2(8!)-2v_2(6!)-2v_2(4!)-2v_2(2)
\]
\[
v_2[H_f(10)]=10v_2(10!)-2v_2(8!)-2v_2(6!)-2v_2(4!)-2v_2(2!)
\]
To conclude:
\[
v_2[H_f(10)]=10v_2(10!)-2\sum_{i=1}^{4}v_2[(2i)!]
\]
So:
\[
v_2[H_f(10)]=80-30=50
\]
And this is our exact result. \\
\section{Theorem 1}
Let $n,p \in\mathbb{N}$
,$p$ a prime integer and $H_f(n)$ is the hyperfactorial function of $n$,\\
hence:

\[
v_p[H_f(n)]=p\Bigl\lfloor\frac{n}{p}\Bigr\rfloor v_p(n!)-p\sum_{i=1}^{\lfloor\frac{n}{p}\rfloor -1} v_p[(pi)!]
\] 
\textit{Proof:}
We've already seen that our formula is working for a particular case. \\
For the proof by induction, let's consider the "n+p" case:
\[
v_p[H_f(n+p)]=v_p[H_f(n)]+(n+1)v_p(n+1)+(n+2)v_p(n+2)+\dots+(n+p)v_p(n+p)\textbf{[1]}
\]
Between "n" and "n+p", obviously one of those numbers is divisible by p.
So only one of this numbers: $v_p(n+1),v_p(n+2),\dots,v_p(n+p)$ \\
is not equal to 0. Let's pay attention on this number, and I define the integer $g \in\mathbb{N}$ such that $0<g\leq p$ e $v_p(n+g)\not=0$
Let's consider this quantity:
\[
(n+g)v_p[(n+g)!]=(n+g)v_p(n+g)+(n+g)v_p(n+g-p)+(n+g)v_p(n+g-2p)+\dots 
\]
Isolating $(n+g)v_p(n+g)$ from the equation above:
\[
(n+g)v_p(n+g)=(n+g)v_p[(n+g)!]-(n+g)v_p(n+g-p)-(n+g)v_p(n+g-2p)-\dots
\]
So:
\[
(n+g)v_p(n+g)=(n+g)v_p[(n+g)!]-(n+g)v_p[(n+g-p)!]
\]
But:
\[
v_p[H_f(n)]=p\Bigl\lfloor\frac{n}{p}\Bigr\rfloor v_p(n!)-p\sum_{i=1}^{\lfloor\frac{n}{p}\rfloor -1} v_p[(pi)!]
\]
And substituting it in [1] the equation just found we'll have:
\[
v_p[H_f(n)]=p\Bigl\lfloor\frac{n}{p}\Bigr\rfloor v_p[(p\Bigl\lfloor\frac{n}{p}\Bigr\rfloor)!]-p\sum_{i=1}^{\lfloor\frac{n}{p}\rfloor -1}v_p[(pi)!]+(n+g)v_p[(n+g)!]-(n+g)v_p[(n+g-p)!]
\]
Notice in this last process that the p-adic valuation is 0 when the argument isn't a multiple of p, hence:
\[
v_p(n!)=v_p[(p\Bigl\lfloor\frac{n}{p}\Bigr\rfloor)!]
\]
But (n+g) and (n+g-p) are divisible by p ,so $0<g\leq p$ and: $(n+g)=p\Bigl\lfloor\frac{n+p}{p}\Bigr\rfloor$
and $(n+g-p)=p\Bigl\lfloor\frac{n}{p}\Bigr\rfloor$ \\
But then:
\[
v_p[H_f(n)]=p\Bigl\lfloor\frac{n}{p}\Bigr\rfloor v_p[(p\Bigl\lfloor\frac{n}{p}\Bigr\rfloor)!]-p\sum_{i=1}^{\lfloor\frac{n}{p}\rfloor -1}v_p[(pi)!]+(p\Bigl\lfloor\frac{n+p}{p}\Bigr\rfloor)v_p[(p\Bigl\lfloor\frac{n+p}{p}\Bigr\rfloor)!]-(p\Bigl\lfloor\frac{n+p}{p}\Bigr\rfloor)v_p[(p\Bigl\lfloor\frac{n}{p}\Bigr\rfloor)!]
\]
Notice that:
\[
p\Bigl\lfloor\frac{n}{p}\Bigr\rfloor v_p[(p\Bigl\lfloor\frac{n}{p}\Bigr\rfloor)!]-p\Bigl\lfloor\frac{n+p}{p}\Bigr\rfloor)v_p[(p\Bigl\lfloor\frac{n}{p}\Bigr\rfloor)!]=-pv_p[(p\Bigl\lfloor\frac{n}{p}\Bigr\rfloor)!]
\]
And this is the $\lfloor\frac{n}{p}\rfloor$ term of the sum. So re-writing this expression we'll obtain:
\[
v_p[H_f(n+p)]=p\Bigl\lfloor\frac{n+p}{p}\Bigr\rfloor v_p[(p\Bigl\lfloor\frac{n+p}{p}\Bigr\rfloor)!]-p\sum_{i=1}^{\lfloor\frac{n}{p}\rfloor }v_p[(pi)!]
\]
And that's the thesis of our induction argument. \\
\\
\\
\section{Expansion of De-Polignac formula}
Let $n,p \in\mathbb{N}$,
$p$ a prime integer, and $H_f(n)$ is the hyperfactorial function \\
of $n$, hence: \\

\[
v_p[H_f(n)]=p\Bigl\lfloor\frac{n}{p}\Bigr\rfloor \sum_{k=1}^{\infty} \Bigl\lfloor{\frac{n}{p^{k}}}\Bigr\rfloor-p\sum_{i=1}^{\lfloor\frac{n}{p}\rfloor -1} \sum_{k=1}^{\infty} \Bigl\lfloor\frac{i}{p^{k-1}}\Bigr\rfloor
\]
We can link our formula with the De-Polignac one \cite{tdn}, so we'll have that:
\[
v_p(n!)=\sum_{k=1}^{\infty} \Bigl\lfloor\frac{n}{p^k}\Bigr\rfloor
\]
But substituting it in the formula proved in the section before we'll have:
\[
v_p[H_f(n)]=p\Bigl\lfloor\frac{n}{p}\Bigr\rfloor \sum_{k=1}^{\infty} \Bigl\lfloor{\frac{n}{p^{k}}}\Bigr\rfloor-p\sum_{i=1}^{\lfloor\frac{n}{p}\rfloor -1} \sum_{k=1}^{\infty} \Bigl\lfloor\frac{pi}{p^{k}}\Bigr\rfloor
\]
And to conclude:
\[
v_p[H_f(n)]=p\Bigl\lfloor\frac{n}{p}\Bigr\rfloor \sum_{k=1}^{\infty} \Bigl\lfloor{\frac{n}{p^{k}}}\Bigr\rfloor-p\sum_{i=1}^{\lfloor\frac{n}{p}\rfloor -1} \sum_{k=1}^{\infty} \Bigl\lfloor\frac{i}{p^{k-1}}\Bigr\rfloor
\]
\section{Asymptotic analysis}
Now we can determine the asymptotic behaviour of our formula as $n\to\infty$ 
\subsection{Lemma 1}
It is known that:
\[
v_p(n!)\sim\frac{n}{p-1}+O(log_p(n))
\]
\subsection{Theorem 2}
\[
v_p[H_f(n)]\sim\frac{n(n+p)}{2(p-1)}
\]
as $n\to\infty$. \\
\textit{Proof:} 
We have that:
\[
v_p[H_f(n)]=p\Bigl\lfloor\frac{n}{p}\Bigr\rfloor v_p(n!)-p\sum_{i=1}^{\lfloor\frac{n}{p}\rfloor -1} v_p[(pi)!]
\]
But $p\lfloor\frac{n}{p}\rfloor\sim n$ and $v_p(n!)\sim\frac{n}{p-1}$ for $n\to\infty$ \\
Furthermore:
\[
\sum_{i=1}^{\lfloor\frac{n}{p}\rfloor -1} v_p[(pi)!]=[v_p(p!)+v_p[(2p)!]+v_p[(3p)!]+\dots+v_p[(p(\lfloor\frac{n}{p}\rfloor -1))!]
\]
Hence:
\[
\sum_{i=1}^{\lfloor\frac{n}{p}\rfloor -1} v_p[(pi)!]\sim \Bigl[\frac{p}{p-1}+\frac{2p}{p-1}+\dots+\frac{\lfloor\frac{n}{p}\rfloor -1}{p-1}\Bigr]=\frac{p}{p-1}\sum_{i=1}^{\lfloor\frac{n}{p}\rfloor -1}i
\]
But:
\[
\frac{p}{p-1}\sum_{i=1}^{\lfloor\frac{n}{p}\rfloor -1}i=\frac{p}{p-1}\frac{(\lfloor\frac{n}{p}\rfloor)(\lfloor\frac{n}{p}\rfloor -1)}{2}\sim \frac{n^2-np}{2p(p-1)}
\]
So:
\[
v_p[H_f(n)]\sim \frac{n^2}{p-1}-p\frac{n^2-np}{2p(p-1)}=\frac{n(n+p)}{2(p-1)}
\]
For example, let's consider $v_7[H_f(1000)]=82390$ \\
With our formula we'll have:
\[
\frac{1000(1000+7)}{2(7-1)}\approx 83917
\]
And this is a nice approximation of the exact result. \\
Notice that:
\[
\lim_{n\to\infty} \frac{v_p(n!)}{v_p[H_f(n)]}=0
\]
In fact:
\[
\lim_{n\to\infty} \frac{v_p(n!)}{v_p[H_f(n)]}=\lim_{n\to\infty} \frac{n}{p-1}\frac{2(p-1)}{n(n+p)}=\lim_{n\to\infty} \frac{2}{n+p}=0 
\]
And noticing that $p>1$. \\
We can do another important observation, let's consider the following limit:
\[
\lim_{n\to\infty} \frac{[v_p(n!)]^2}{v_p[H_f(n)]}=l
\]
But:
\[
[v_p(n!)]^2\sim\frac{n^2}{(p-1)^2}
\]
Hence:
\[
\lim_{n\to\infty} \frac{[v_p(n!)]^2}{v_p[H_f(n)]}=\lim_{n\to\infty} \frac{n^2}{(p-1)^2}\frac{2(p-1)}{n(n+p)}=\frac{2}{p-1}
\]
But then we'll have:
\[
v_p(n!)\sim \sqrt{\frac{2}{p-1}v_p[H_f(n)]}
\]
And:
\[
v_p[H_f(n)]\sim\frac{(p-1)[v_p(n!)]^2}{2}
\]
\section{Graphs}
\begin{figure}[H]
    \centering
    \includegraphics[width=6cm]{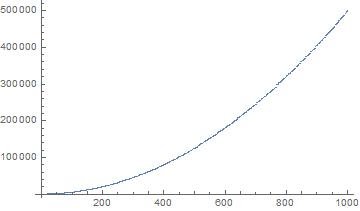}
    \caption{Graph of $v_2[H_f(n)]$ from $n$=1 to $n$=1000}
    \label{fig:1}
\end{figure}
\begin{figure}[H]
    \centering
    \includegraphics[width=6cm]{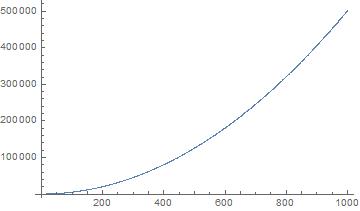}
    \caption{Graph of $\frac{n^{2}+2n}{2}$ from $n$=1 to $n$=1000}
    \label{fig:1}
\end{figure}
\begin{figure}[H]
    \centering
    \includegraphics[width=6cm]{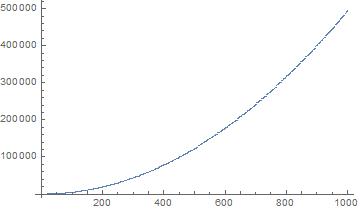}
    \caption{Graph of $\frac{[v_2(n!)]^{2}}{2}$ from $n$=1 to $n$=1000}
    \label{fig:1}
\end{figure}
\bibliographystyle{plain}
\bibliography{Bibliografia.bib}
\end{document}